\documentclass[12pt]{article}
\usepackage{amssymb,latexsym}
\oddsidemargin=\evensidemargin
\advance\topmargin by-\headheight
\def\A{{\cal A}}
\def\AGL{\mathop{\rm AGL}\nolimits}
\def\AGaL{\mathop{\rm A\Gamma L}\nolimits}
\def\aut{\mathop{\rm Aut}\nolimits}

\def\dim{\mathop{\rm dim}\nolimits}

\def\End{\mathop{\rm End}\nolimits}
\def\F{\mathbb F}

\def\H{{\cal H}}

\def\ker{\mathop{\rm ker}\nolimits}

\def\lg{\langle}

\def\O{{\cal O}}
\def\orb{\mathop{\rm Orb}\nolimits}

\def\rad{\mathop{\rm rad}\nolimits}
\def\rg{\rangle}
\def\rk{\mathop{\rm rk}\nolimits}
\def\S{{\cal S}}

\def\sym{\mathop{\rm Sym}\nolimits}

\def\wh{\widehat}

\def\ZZ{{\mathbb Z}}

\def\proof{{\bf Proof}.\ }
\def\bull{\vrule height .9ex width .8ex depth -.1ex }

\def\SuSS{\stepcounter{subsection}{\bf\thesubsection{.}}
\addtocounter{subsection}{-1}\refstepcounter{subsection}}
\def\eSuSS{\vspace{2mm}\SuSS}

\newtheorem{formula}{}[section]
\newtheorem{proposition}[formula]{Proposition}
\newtheorem{definition}[formula]{Definition}
\newtheorem{corollary}[formula]{Corollary}
\newtheorem{remark}[formula]{Remark}
\newtheorem{lemma}[formula]{Lemma}
\newtheorem{theorem}[formula]{Theorem}

\def\thrm{\begin{theorem}}
\def\thrml#1{\begin{theorem}\label{#1}}
\def\ethrm{\end{theorem}}
\def\prpstn{\begin{proposition}}
\def\prpstnl#1{\begin{proposition}\label{#1}}
\def\eprpstn{\end{proposition}}
\def\rmrk{\begin{remark}}
\def\rmrkl#1{\begin{remark}\label{#1}}
\def\ermrk{\end{remark}}
\def\dfntn{\begin{definition}}
\def\dfntnl#1{\begin{definition}\label{#1}}
\def\edfntn{\end{definition}}
\def\nmrt{\begin{enumerate}}
\def\enmrt{\end{enumerate}}
\def\tm#1{\item[{\rm (#1)}]}
\def\qtn{\begin{equation}}
\def\qtnl#1{\begin{equation}\label{#1}}
\def\eqtn{\end{equation}}
\def\lmm{\begin{lemma}}
\def\lmml#1{\begin{lemma}\label{#1}}
\def\elmm{\end{lemma}}
\def\crllr{\begin{corollary}}
\def\crllrl#1{\begin{corollary}\label{#1}}
\def\ecrllr{\end{corollary}}

\begin{document}
\title{A new look at the Burnside-Schur theorem}
\author{
Sergei Evdokimov \\[-1pt]
\small St. Petersburg Institute \\[-3pt]
\small for Informatics and Automation RAS\\[-3pt]
{\tt \small evdokim@pdmi.ras.ru }
\thanks{Partially supported by RFFI, grants 01-01-00219, 03-01-00349}
\and
Ilia Ponomarenko\\[-1pt]
\small Steklov Institute of Mathematics\\[-3pt]
\small at St. Petersburg \\[-3pt]
{\tt \small inp@pdmi.ras.ru}
\thanks{Partially supported by RFFI, grants, 01-01-00219, 03-01-00349,
NSH-2251.2003.1}
}
\date{23.10.2003}
\maketitle

\begin{abstract}
The famous Burnside-Schur theorem states that every primitive finite
permutation group containing a regular cyclic subgroup is either 2-transitive
or isomorphic to a subgroup of a 1-dimensional affine group of prime degree.
It is known that this theorem can be expressed as a statement on Schur
rings over a finite cyclic group. Generalizing the latters we introduce
Schur rings over a finite commutative ring and prove an analog of this
statement for them. Besides, the finite local commutative rings are
characterized in the permutation group terms.
\end{abstract}

\section{Introduction}

\SuSS\label{l021003}
The starting point of the paper is the following statement.
\vspace{2mm}

{\bf Theorem (Burnside-Schur).} Every primitive finite permutation group
containing a regular cyclic subgroup is either 2-transitive or
permutationally isomorphic to a subgroup of the affine group $\AGL_1(p)$ where
$p$ is a prime.\footnote{In the prime degree case the assumption on the
existence of a regular cyclic subgroup is unnecessary.}\bull
\vspace{2mm}

In fact, Burnside~\cite[pp.339-343]{B} proved the theorem for the
permutation groups of prime power degree and conjectured that the
statement would still true when the regular subgroup in question was
abelian, not only cyclic. Schur \cite{S33} disproved this conjecture
and showed that every primitive permutation group of composite degree
containing a regular cyclic subgroup is 2-transitive. This result was
generalized by Wielandt in \cite[Theorem~25.4]{W64} (see also
Corollary~\ref{crw} below). 

In contrast to the method of Burnside who employed the character theory
the idea of Schur was as follows. Let $\Gamma$ be a permutation group
containing a regular subgroup~$G$. Then the submodule $\A(\Gamma)$ of
the group ring of~$G$ spanned by the orbits of a one-point stabilizer
of~$\Gamma$ is a subring of the latter ring (after the choice of the
point these orbits are treated in a natural way as subsets of~$G$).
The ring $\A(\Gamma)$ is a special case of a Schur ring (or S-ring) over
the group~$G$ (for the definition
of S-rings and basic facts about them see 
Section~\ref{sec2}).\footnote{ The term was proposed by 
Wielandt in~\cite{W64} where the theory of S-rings was developed.} Schur
observed that the group~$\Gamma$ is 2-transitive (resp. primitive)
iff the S-ring $\A(\Gamma)$ is of rank~2 (resp. primitive). Thus the
Burnside-Schur theorem when the degree of~$\Gamma$ is a composite
number is an immediate consequence of the following theorem
proved in fact by Schur: every primitive S-ring over a cyclic group of
composite order is of rank~2.

An essential ingredient of Schur's proof is the following theorem on
multipliers: every S-ring over a finite abelian group~$G$ is invariant with
respect to the group~$K_G$ consisting of the automorphisms of~$G$ induced
by raising to powers coprime to~$|G|$. If the group~$G$ is
cyclic, then it can be treated as
the additive group of the ring $\ZZ_n$ of integers modulo~$n$ where
$n=|G|$. In this case $K_G$ is identified with the
multiplications by the units of~$\ZZ_n$. Thus in the cyclic case the
multiplier theorem states that every S-ring over~$G$ is an S-ring
{\it over the ring} $\ZZ_n$ in the sense of the following definition.
\vspace{2mm}

{\bf Definition.} {\it Let $\A$ be an S-ring over a group~$G$ and $R$
a finite commutative ring. We say that~$\A$ is an S-ring over the ring
$R$ if $G=R^+$ and $\A$ is invariant with respect to the group
$K_R\le\aut(G)$ induced by the action of~$R^\times$ on $R^+$ by
multiplications.}
\vspace{2mm}

{\bf Example.} Let $R$ be a finite commutative ring and $K\le K_R$. Then,
obviously, the set of all $K$-invariant elements of $\ZZ[R^+]$ 
forms an S-ring over the ring~$R$. It is
called a {\it cyclotomic} one. The basic sets of this ring are exactly
the orbits of~$K$ on~$R$. Since the group $K_R$ is regular on~$R^\times$,
the cyclotomic rings are in 1-1 correspondence to the subgroups of~$R^\times$.
The translation association schemes corresponding to cyclotomic S-rings were 
studied in~\cite{GC}. Let now $R$ is a field of order~$q$ and 
$\Gamma\le\sym(R)$ a permutation group containing the translations
of~$R$. Suppose that $\A(\Gamma)$ is a cyclotomic ring of rank greater
than~2. Then $\Gamma$ is a subgroup of the automorphism group of the
translation association scheme (Cayley ring) corresponding to~$\A(\Gamma)$
(see e.g.~\cite{EP}). So $\Gamma\le\AGaL_1(q)$ by~\cite[p.389]{BCN}.

It should be remarked that
the group $K_G$ as in the Schur theorem on multipliers coincides with the
center of $\aut(G)$. If $G$ is cyclic, then obviously $K_G=\aut(G)$
and so any S-ring over~$G$ admits the maximal possible multiplier group.
This probably explains the fact that there is a nice theory of S-rings
over finite cyclic groups~\cite{LM,EP}. On the contrary, $K_G<\aut(G)$ for 
any noncyclic~$G$ (in this case $K_G$ is not even a maximal abelian subgroup 
of $\aut(G)$) and the theory of S-rings over arbitrary abelian groups is far 
from being completed. The present
paper is the first step to reconstruct the main features of the cyclic case
theory for the S-rings admitting an appropriate sufficiently large
multiplier group containing $K_G$ (e.g. $K_R$ where $R$ is a finite
commutative ring with $R^+=G$). In subsequent papers we plan to study
S-rings over finite commutative rings  and describe them completely at
least for the products of Galois rings of pairwise coprime characteristics.

\eSuSS
To formulate the main results of the paper we recall some basic facts on
finite rings (see e.g.~\cite{MD}). Let $R$ be a finite ring with identity.
Then
$$
R=\prod_pR_p
$$
where $p$ runs over all prime divisors of~$|R|$ and $R_p$ is
a {\it primary component} of~$R$, i.e. a subring of~$R$ such that
$(R_p)^+$ is the Sylow $p$-subgroup of~$R^+$. Moreover, each commutative
primary component
of~$R$ is a direct product of local rings (i.e. ones the non-units of which
form an ideal). Below an S-ring over a finite commutative ring~$R$ is
called {\it quasiprimitive} if it is $K_R$-primitive (see the definition
before Theorem~\ref{thrm2}). If $R$ is not isomorphic to the product of
two rings one of which is $\ZZ_2\times\ZZ_2$, then it is generated by
the units (see~\cite[p.406]{MD}) and so in this case an S-ring $\A$ over~$R$
is quasiprimitive iff $\{0\}$ and~$R$ are the only ideals of~$R$ that
are also $\A$-subgroups of~$R^+$. Now the following statement can be 
considered as a generalization of the Burnside-Schur theorem.

\thrml{thrm1}
Let $R$ be a finite commutative ring with identity. If every primary
component of~$R$ is a local ring, then each quasiprimitive S-ring
over~$R$ is either of rank~$2$ or cyclotomic. In the latter case
$R$ is a field.
\ethrm

To see that Theorem~\ref{thrm1} really generalizes the Burnside-Schur theorem
let $\Gamma$ be a primitive permutation group containing a regular cyclic
subgroup~$G$. Then by the Schur theorem on multipliers one can assume
(without loss of gnerality) that $\A(\Gamma)$ is an S-ring over the
ring~$\ZZ_n$ where $n=|G|$. Since $\Gamma$ is a primitive group, $\A(\Gamma)$
is a primitive S-ring over the group $\ZZ_n^+$, which means that it is quasiprimitive
S-ring over the ring~$\ZZ_n$. Besides, every primary component of~$\ZZ_n$ is obviously a
local ring. Thus the statement of Theorem~\ref{thrm1} (with $R=\ZZ_n$) holds
for the S-ring $\A(\Gamma)$. If $\rk(\A(\Gamma))=2$, then $\Gamma$ is
2-transitive. Otherwise, $\A(\Gamma)$ is cyclotomic and $n=p$ is a prime.
So $\Gamma\le\AGL_1(p)$ by~\cite[Proposition~12.7.5]{BCN} (see the
example in Subsection~\ref{l021003}).

Theorem~\ref{thrm1} is an immediate consequence of Theorem~\ref{thrm3} below
which will be deduced from the following theorem on groups to be proved
in Section~\ref{s300703}. Below an S-ring $\A$ over a group~$G$ is called
{\it $K$-primitive} where $K\le\aut(G)$, if $\{1\}$ and~$G$ are the only
$K$-invariant $\A$-subgroups of~$G$ (for $K=\{1\}$ this means that $\A$ is 
primitive).

\thrml{thrm2}
Let $G$ be a finite abelian group and $K\le\aut(G)$. Suppose that there
exists a Sylow $p$-subgroup $P$ of $G$ such that
\nmrt
\tm{1} $K\cap\aut(P)$ is an abelian group,
\tm{2} $P$ is the disjoint union of an orbit of $K\cap\aut(P)$ and a
$K$-invariant subgroup of~$P$.
\enmrt
Then each $K$-primitive
$K$-invariant S-ring over $G$ is either of rank~$2$ or Cayley
isomorphic to a cyclotomic ring over a field with $K$ going to
the multiplications by nonzero elements of the field.
\ethrm

\crllrl{crw}{\rm\bf (Wielandt).}
Each primitive permutation group containing an abelian regular subgroup
of composite order which has a cyclic Sylow subgroup, is 2-transitive. 
\ecrllr
\proof Let $G$ be an abelian group of composite order such that one of its
Sylow subgroups, say $P$, is cyclic. It suffices to prove that every
primitive S-ring over~$G$ is of rank~2. Set~$K=K_G$.  
Then by the Schur theorem on multipliers any S-ring
over~$G$ is $K$-invariant. Moreover, it is easily seen that $P$ is the disjoint
union of the $K$-orbit containing a generator of $P$, and the subgroup of~$P$ 
of prime index. So by Theorem~\ref{thrm2} every $K$-primitive S-ring over~$G$
is either of rank~$2$ or cyclotomic over a field. Since $|G|$ is a
composite number, the latter case is impossible. To complete the proof
it suffices to note that obviously an S-ring over~$G$ is $K$-primitive
iff it is primitive.\bull

Let $R$ be a finite commutative ring with identity at least one primary
component of which, say $R_p$, is a local ring. Then the hypothesis of
Theorem~\ref{thrm2} is satisfied for $G=R^+$, $K=K_R$ and $P=(R_p)^+$
(indeed, $K\cap\aut(P)=K_{R_p}$, $(1_{R_p})^{K_{R_p}}=(R_p)^\times$ and
$P=(R_p)^\times\cup\rad(R_p)$). Besides, in accordance with our definitions the
quasiprimitive S-rings over~$R$ are exactly the $K$-primitive $K$-invariant
S-rings over~$G$. Thus the following statement is a specialization of
Theorem~\ref{thrm2}

\thrml{thrm3}
Let $R$ be a finite commutative ring with identity. If at least one primary
component of~$R$ is a local ring, then each quasiprimitive S-ring
over~$R$ is either of rank~$2$ or cyclotomic. In the latter 
case $R$ is a field.\bull
\ethrm

As the following example shows the locality of some primary
component of the ring~$R$ in Theorem~\ref{thrm3} is essential. Indeed,
let $R=\ZZ_p\times\ZZ_p$ where $p>2$ is a prime. Set
$$
X_0=\{(0,0)\},\quad 
X_1=(\ZZ_p^\times\times\{0\})\cup(\{0\}\times\ZZ_p^\times),\quad
X_2=\ZZ_p^\times\times\ZZ_p^\times.
$$
Then, obviously, $X_0$, $X_1$, $X_2$ are the basic sets of an S-ring~$\A$
over the group~$R^+$. Since $R^\times=\ZZ_p^\times\times\ZZ_p^\times$
and $\ZZ_p\times\{0\}$ and $\{0\}\times\ZZ_p$ are the only proper ideals 
of~$R$, one can see that $\A$ is a quasiprimitive S-ring
over the ring~$R$. It remains to note that $\rk(\A)=3$ and the ring~$R$
is not local.

Let $G=P$ be a finite abelian $p$-group. As we saw above the hypothesis of 
Theorem~\ref{thrm2} is satisfied when $K=K_R$ where $R$ is
a local commutative ring on~$P$, i.e. one with $R^+=P$. In fact this
is the only possible case because as $R$ runs over all such rings,
$R^\times$ runs over all groups~$K$ as in this theorem (see 
Theorem~\ref{thrm4}). To be more precise, let $K\le\aut(P)$ be an
abelian group and $e\in P$. We say that $(K,e)$ is a {\it local pair}
on~$P$ if the set $P\setminus O$ is a subgroup of~$P$ where
$O=e^K$. It is easily seen that this subgroup is uniquely determined
by~$K$ and the group~$K$ acts regular and faithfully on~$O$.

\thrml{thrm4}
Let $P$ be a finite abelian $p$-group. Then the mapping
\qtnl{l240703}
R\mapsto (K_R,1_R)
\eqtn
establishes a 1-1 correspondence between the set of all local commutative
rings on~$P$ and the set of all local pairs on~$P$. Moreover, two such rings
are isomorphic iff the corresponding subgroups are conjugate in~$\aut(P)$.
\ethrm

The proof of Theorem~\ref{thrm4} as well as Theorem~\ref{thrm2} is contained
in Section~\ref{s300703}. In fact, the latter theorem is deduced from the
first one and the theorem on a separating subgroup proved in
Section~\ref{s130903}. Section~\ref{sec2} contains necessary definitions and
facts on S-rings. All undefined terms and results concerning permutation groups 
can be found in~\cite{DM}.

\eSuSS
We complete the introduction by making some remarks on a topic closely
related to the contents of the paper. Following~\cite{W64} a finite 
group~$G$ is called a {\it B-group} if every primitive group containing
a regular subgroup isomorphic to~$G$ is 2-transitive. In fact, the
Burnside-Schur theorem states in particular that a cyclic group of
composite order is a B-group (see also Corollary~\ref{crw}). As P.~Cameron
observed in~\cite{C} the classification of finite simple groups implies
that for almost all positive integers~$n$, every group of order~$n$ is
a B-group. At the same time the problem of the classification of B-groups
(even abelian ones) is still open. In this connection it should be
remarked that most of B-groups~$G$  listed in~\cite{W64} including
those as in Corollary~\ref{crw}
satisfy {\it a priori} a more strong condition: every primitive S-ring
over~$G$ is of rank~2 (in fact there exist S-rings even over a cyclic
group that do not come from permutation groups, see~\cite{EP01}). Obviously,
no group of odd prime order satisfies this condition. However, the
following is true.

\prpstnl{p23}
For any abelian group~$G$ there exists an abelian group
$K\le\aut(G)$ such that each $K$-primitive $K$-invariant S-ring over~$G$
is of rank~2.
\eprpstn
\proof By Theorem~\ref{thrm1} 
it suffices to verify that each abelian $p$-group~$P$ of composite order
is isomorphic to the additive group of a local commutative ring which is
not a field. By Theorem~\ref{thrm4} all we need is to present a local
pair $(K,e)$ on~$P$ such that $P\setminus e^K\ne \{0\}$. Let $p^m$
be the exponent of~$P$ and $P=\ZZ^+_{p^m}\times P'$ for some~$P'$. Then
treating $P$ as a $\ZZ_{p^m}$-module one can take $e=(1,0)$ and
$K=\{x\mapsto ax+x_1b,\ x\in P:\ a\in\ZZ^\times_{p^m},\ b\in P'\}$
where $x_1$ is the first coordinate of~$x$. It is easy to see that
the group~$K$ is abelian.\bull

Proposition \ref{p23} shows that the study of pairs $(G,K)$ as in it
is seemingly more fruitful than that of B-groups. For instance, it would be 
interesting  to identify B-rings, i.e. finite commutative rings $R$ for 
which $(R^+,K_R)$ is such a pair (a special class of B-rings can be derived
from Theorem~\ref{thrm3}).

{\bf Notation.} As usual by $\ZZ$ we denote the ring of integers.

For a ring~$R$ with identity we denote by $R^+$, $R^\times$ and
$\rad(R)$ the additive and
multiplicative groups of~$R$ and the radical of~$R$ respectively.

The group of all permutations of a set $V$ is denoted by $\sym(V)$. If
$\Gamma\le\sym(V)$, then $\orb(\Gamma,V)$ denotes the set of all
orbits of the group~$\Gamma$ on~$V$.

\section{S-rings}
\label{sec2}

\SuSS
Let $G$ be a finite group. A subring~$\A$ of the group ring~$\ZZ[G]$ is
called a {\it Schur ring} (briefly {\it S-ring}) over~$G$ if it has a
(uniquely determined) $\ZZ$-base consisting of the elements
$\xi(X)=\sum_{x\in X}x$ where $X$ runs over a family $\S=\S(\A)$ of pairwise
disjoint nonempty subsets of~$G$ such that
$$
\{1\}\in\S,\quad
\bigcup_{X\in\S}X=G\quad
\textstyle{\rm and}\quad
X\in\S\ \Rightarrow\ X^{-1}\in\S.
$$
We call the elements of $\S$ {\it basic} sets of~$\A$ and denote by $\S^*(\A)$
the set of all unions of them and by~$\H(\A)$ the set of all {\it $\A$-subgroups}
of~$G$, i.e. those belonging to~$\S^*(\A)$. For $K\le\aut(G)$ we set
$\H_K(\A)=\{H\in\H(\A):\ H^K=H\}$. The basic set of~$\A$ that
contains $x\in G$ is denoted by $[x]$. The number $\rk(\A)=\dim_\ZZ(\A)$ 
is called
the {\it rank} of~$\A$. Two S-rings $\A$ over~$G$ and~$\A'$ over~$G'$ are
called {\it Cayley isomorphic} if there exists a group isomorphism $f:G\to G'$
such that $\A'$ equals the image of~$\A$ with respect to the isomorphism
from $\ZZ[G]$ to $\ZZ[G']$ induced by~$f$.

The first statement of the following lemma was in fact proved in a different
form in~\cite[Theorem~23.9]{W64} which also contains the proof of the
second statement for $K=\{1\}$.

\lmml{t100703b}
Let $G$ be a finite abelian group, $\A$ an S-ring over $G$ and
$p$ a prime dividing $|G|$. Then for any $X\in\S(\A)$ the
following statements hold:
\nmrt
\tm{1} $X^{[p]}\in\S^*(\A)$ where 
$X^{[p]}=\{x^p:\ x\in X,\ |xE\cap X|\not\equiv 0\pmod p\}$ with
$E=\{g\in G:\ g^p=1\}$. 
\tm{2} if the ring $\A$ is $K$-invariant and $K$-primitive for
some group $K\le\aut(G)$, then $X^{[p]}=\{1\}$; in particular,
$|xE\cap X|\equiv 0\pmod p$ for all $x\in G\setminus E$.
\enmrt
\elmm
\proof Let $X\in\S(\A)$. Then using a well-known property of binomial
coefficients we have
\qtnl{l010803}
\xi(X)^p=(\sum_{x\in X}x)^p\equiv\sum_{x\in X}x^p\pmod p.
\eqtn
On the other hand, given a coset $C\in G/E$, the
element $g^p$ does not depend on the choice of $g\in C$. Denote it
by~$h_C$. We observe that the mapping $C\mapsto h_C$ is a bijection
from $G/E$ onto~$G^p$. So
\qtnl{l010803a}
\sum_{x\in X}x^p=\sum_{C\in G/E}|C\cap X|\cdot h_C\equiv
\sum_{C\in G/E,\atop |C\cap X|\not\equiv 0\pmod p}|C\cap X|\cdot h_C
\pmod p.
\eqtn
By \cite[Proposition~22.3]{W64} formulas~(\ref{l010803}) and~(\ref{l010803a})
imply that the set $\{h_C:\ C\in G/E,\ |C\cap X|\not\equiv 0\pmod p\}$
belongs to $\S^*(\A)$. Since this set equals $X^{[p]}$, statement~(1)
is proved. To prove statement~(2) suppose that
the ring $\A$ is $K$-invariant and 
$K$-primitive for some $K\le\aut(G)$. It suffices to verify that
$X^{[p]}=\{1\}$. To do this denote by~$H$ the subgroup of~$G$ generated
by the sets $(X^{[p]})^k$, $k\in K$. Then $H\in\H_K(\A)$ by statement~(1)
and the $K$-invariance of~$\A$. On the other hand, $H\le G^p<G$. So
$H=\{1\}$ due to the $K$-primitivity of~$\A$. Since $X^{[p]}\subset H$,
we are done.\bull

\eSuSS
Let $G$ be a finite abelian group and $\wh G$ the group dual to~$G$.
For $\sigma\in\aut(G)$ and $\chi\in\wh G$ we define a function
$\chi^{\wh\sigma}$ on $\wh G$ by $\chi^{\wh\sigma}(g)=\chi(g^{\sigma^{-1}})$,
$g\in G$. It is easily seen that $\chi^{\wh\sigma}\in\wh G$. Moreover,
the mapping $\chi\mapsto\chi^{\wh\sigma}$, $\chi\in\wh G$, belongs to
$\aut(\wh G)$ and the mapping
$$
\aut(G)\to\aut(\wh G),\quad\sigma\mapsto\wh\sigma,
$$
is a group isomorphism. This enables us to identify $\aut(G)$ with
$\aut(\wh G)$ and treat any subgroup~$K$ of $\aut(G)$ as a subgroup
of $\aut(\wh G)$. Clearly, $H$ is a $K$-invariant subgroup of~$G$ iff
$H^\bot$ is a $K$-invariant subgroup of~$\wh G$ where
$H^\bot=\{\chi\in\wh G:\ H\le\ker(\chi)\}$.

Let $\A$ be an S-ring over the group~$G$. Denote by $\wh S$ the set of all
classes of the equivalence relation on~$\wh G$ defined as follows:
$\chi_1\sim\chi_2$
iff the extensions of~$\chi_1$ and $\chi_2$ to $\ZZ[G]$ coincide on~$\A$.
Then the $\ZZ$-submodule of $\ZZ[\wh G]$ spanned by the
elements $\xi(X)$, $X\in\wh S$, is an S-ring over~$\wh G$ 
(see~\cite[Theorem~6.3]{BI}). This ring is called {\it dual} to~$\A$
and denoted by~$\wh\A$. Obviously, $\S(\wh\A)=\wh S$. Moreover,
$\rk(\A)=\rk(\wh\A)$.

\thrml{t110703a}
Let $\A$ be an S-ring over a finite abelian group $G$ and $K\le\aut(G)$.
Then
\nmrt
\tm{1} the ring $\A$ is $K$-invariant iff the dual ring $\wh\A$ is
$K$-invariant.
\tm{2} $H\in\H_K(\A)$ iff $H^\bot\in\H_K(\wh\A)$; in particular,
the ring $\A$ is $K$-primitive iff the dual ring $\wh\A$ is
$K$-primitive.
\enmrt
\ethrm
\proof If the ring $\A$ is $K$-invariant, then given $\sigma\in K$ we
have:  $\chi_1\sim\chi_2$ iff $\chi_1^\sigma\sim\chi_2^\sigma$. So the
ring~$\wh\A$ is $K$-invariant. The converse statement follows from
the equality $\wh{\wh\A}=\A$. To prove the first part of statement~(2)
we observe that given $H\le G$ and $\chi\in\wh G$ we have:
$H\le\ker(\chi)$ iff $\sum_{h\in H}\chi(h)=|H|$ (we used that
$|\chi(h)|=1$ for all~$h$ and $\chi(1)=1$). Let now $H\in\H_K(\A)$. 
Then the above observation implies that
$\chi\in H^\bot$ iff $\chi'\in H^\bot$ whenever $\chi\sim\chi'$. So 
$[\chi]\subset H^\bot$ for all $\chi\in H^\bot$, and hence 
$H^\bot\in\H_K(\wh\A)$. The second part of statement~(2) is the consequence
of the first one and the obvious equalities $\{1\}^\bot=\wh G$ and
$G^\bot=\{1\}$.\bull 

\section{Theorem on a separating subgroup}
\label{s130903}

In this section we prove a statement on S-rings over a finite group
which is a key one for the proof of Theorem~\ref{thrm2}. In this
connection it is worth remarking that the proofs of Wielandt's theorem
(Corollary~\ref{crw}) from~\cite{BCN} and~\cite{DM} go back to its
original proof (see~\cite[Theorem~25.4]{W64}). A detailed
analysis shows that in fact most part of Wielandt's proof deals with
a special case of Theorem~\ref{t100703} below. In its turn this theorem 
is a consequence of a more general result on association schemes the 
proof of which is outside the scope of the present paper.

The following definition is taken from~\cite{EP}.
Given a nonempty subset $X$ of a finite group~$G$ the group
$$
\rad(X)=\{g\in G:\ gX=Xg=X\}
$$
is called the {\it radical} of~$X$.\footnote{Not to mix it up with the 
radical of a ring.} It is the largest subgroup of~$G$ such that~$X$ 
is a union of the left as well as right cosets by this subgroup. 
Besides, obviously $\rad(X)\subset\lg X\rg$ where $\lg X\rg$ is a
subgroup of~$G$ generated by~$X$. If $X\in\S^*(\A)$ where $\A$ is an 
S-ring over $G$, then $\rad(X)$ and $\lg X\rg$ are $\A$-subgroups 
of~$G$.

\thrml{t100703}
Let $X$ be a basic set of an S-ring~$\A$ over a finite group $G$. 
Suppose that
$$
\lg X\cap H\rg\le\rad(X\setminus H)
$$
for some subgroup $H$ of~$G$ such that $X\cap H\ne\emptyset$ and
$X\setminus H\ne\emptyset$. Then $X=\lg X\rg\setminus\rad(X)$
with $\rad(X)\le H\le\lg X\rg$.
\ethrm
\proof For $A,B,C\subset G$ set $\xi_A=\xi(A)$ and 
$\xi_{A,B,C}=(\xi_A\xi_B)\circ\xi_C$. Since $X\in\S(\A)$,
we have
\qtnl{l050803}
\xi_{X,X^{-1},X}=a_X\xi_X=a_X\xi_Y+a_X\xi_Z
\eqtn
for some integer $a_X\ge 0$, where $Y=X\cap H$ and $Z=X\setminus H$.
Obviously, $\xi_{Y,Y^{-1},Z}=\xi_{Y,Z^{-1},Y}=\xi_{Z,Y^{-1},Y}=0$
and hence
\qtnl{l130703a}
\begin{array}{l}
\xi_{X,X^{-1},X}=
((\xi_Y+\xi_Z)(\xi_{Y^{-1}}+\xi_{Z^{-1}}))\circ(\xi_Y+\xi_Z)=\\
(\xi_{Y,Y^{-1},Y}+\xi_{Z,Z^{-1},Y})+
(\xi_{Y,Z^{-1},Z}+\xi_{Z,Y^{-1},Z}+\xi_{Z,Z^{-1},Z}).
\end{array}
\eqtn
From the hypothesis of the theorem it follows that 
\qtnl{l050803a}
\xi_{Z,Z^{-1},Y}=|Z|\xi_Y,\quad
\xi_{Y,Z^{-1},Z}=\delta|Y|\xi_Z,\quad
\xi_{Z,Y^{-1},Z}=|Y|\xi_Z
\eqtn
$\delta=\delta_{X,X^{-1}}$ is the Kronecker delta. Due to~(\ref{l050803})
and~(\ref{l130703a}) this implies that $\xi_{Y,Y^{-1},Y}=a_Y\xi_Y$ and
$\xi_{Z,Z^{-1},Z}=a_Z\xi_Z$ for some integers $a_Y\ge 0$ and
$a_Z\ge 0$. Thus, 
\qtnl{l130703b}
a_X=|Z|+a_Y=(\delta+1)|Y|+a_Z.
\eqtn 

We observe that $a_Y=|Y\cap gY|$ for $g\in Y$ and $a_Z=|Z\cap gZ|$ for 
$g\in Z$. Now let us prove that $\delta=1$ and
\qtnl{l050803c}
a_Y\ge 2|Y|-|H|,\quad a_Z\le |Z|-|H|.
\eqtn
Since $H\cap Z=\emptyset$, we have $gH\cap gZ=\emptyset$ for $g\in Z$.
So $a_Z\le |Z|-|H|$ because $gH\subset Z$. This implies that 
$\delta=1$, for otherwise
$$
a_X=|Y|+a_Z\le |Y|+|Z|-|H|<|Z|\le|Z|+a_Y=a_X
$$
(we used (\ref{l130703b}) and the inclusion $Y\subset H\setminus\{1\}$).
In particular, $Y=Y^{-1}$ and $Z=Z^{-1}$. Moreover, the latter equality
implies that $a_Z=|Z|-|H|$
iff $Z\cup gZ=Z\cup H$ for all $g\in Z$, i.e. iff $Z\cup H$ is a
subgroup of~$G$. Next, by the inclusion-exclusion principle we have
$$
a_Y=|Y\cap gY|=2|Y|-|Y\cup gY|\ge 2|Y|-|H|,\quad g\in Y,
$$
with the equality attained iff $Y\cap gY=H$ for all $g\in Y$, i.e. iff
$H\setminus Y$ is a subgroup of $H$ (indeed, otherwise 
$H=Y\cup (h_1h_2^{-1})Y$ for some $h_1,h_2\in H\setminus Y$ with
$h_1h_2\in Y$, whence $h_1\in(h_1h_2^{-1})Y$ which is impossible).

From (\ref{l130703b}) and the equality $\delta=1$ it follows that
$$
0=(a_Y-2|Y|+|H|)+(|Z|-|H|-a_Z).
$$
By (\ref{l050803c}) both expressions in the brackets are nonnegative.
So they equal to~$0$. Due to the above paragraph this means
that $H\setminus Y$ and $Z\cup H$ are subgroups of~$G$. Since
obviously $X=(H\setminus Y)\cup(Z\cup H)$, we conclude that
$\lg X\rg=Z\cup H$ and $\rad(X)=H\setminus Y$, which completes
the proof.\bull

If the hypothesis of Theorem~\ref{t100703} is satisfied, we say that
$H$ {\it separates} $X$. In this case, obviously, $H$ is a proper
subgroup of~$G$ and $X\ne\{1\}$. Moreover, from Theorem~\ref{t100703}
it follows that $X$ is uniquely
determined by~$H$ (in fact, $X$ is the set difference of the smallest
$\A$-subgroup of~$G$ containing~$H$ and the largest $\A$-subgroup of~$G$ 
contained in~$H$). Denote by $\H_{sep}(\A)$ the set of all subgroups
of~$G$ each of which separates some basic set of the S-ring~$\A$. 
Obviously, $\H_{sep}(\A)\cap\H(\A)=\emptyset$.

\crllrl{l300703}
Let $G$ be a finite group, $K\le\aut(G)$ and $\A$ a $K$-primitive
$K$-invariant S-ring over $G$. Then $\rk(\A)=2$ whenever
$\H_{sep}(\A)$ contains a $K$-invariant group.
\ecrllr
\proof Immediately follows from Theorem~\ref{t100703} and the lemma
below to be also used in Section~\ref{s300703}.

\lmml{t100703a}
In the conditions of Corollary~\ref{l300703} given a proper $K$-invariant 
subgroup $H$ of $G$ we have:
\qtnl{l150703c}
H'\le H\ \Longrightarrow\ H'=\{1\},\quad
H'\ge H\ \Longrightarrow\ H'=G
\eqtn
for any group $H'\in\H(\A)$.
\elmm
\proof Let $H'\in\H(\A)$.
Set $H''$ to be the group generated by the  groups $(H')^k$ for $k\in K$
if $H'\le H$, and  to be the intersection of the same groups if $H'\ge H$.
Then $H''\in\H_K(\A)$ with $H''\ne G$ in the first case and $H''\ne\{1\}$
in the second. Thus we are done by the $K$-primitivity of~$\A$.\bull

\section{Proofs of theorems}
\label{s300703}

\SuSS {\bf Proof of Theorem~\ref{thrm4}.}
Let us prove the injectivity of the mapping~(\ref{l240703}). It suffices
to verify that the multiplication in a local commutative ring $R$ on~$P$
is uniquely determined by the corresponding local pair $(K,e)=(K_R,1_R)$.
Let $r\in R$. If $r\in R^\times$, then
$x\cdot r=x^{k(r)}$ for all $x$ where $k(r)$ is the element of~$K$ induced
by~$r$.  We observe that $k(r)$ depends only on~$r$ but not on the
multiplication law in the ring~$R$. (Indeed, $e^{k(r)}=r$ because $e=1_R$,
and due to the regularity there is the only
element of~$K$ taking~$e$ to~$r$.) If $r\not\in R^\times$, then
due to the locality of~$R$, the element $e+r$ belongs to~$R^\times$
and $x\cdot r=x\cdot (e+r)-x$. Thus the injectivity is proved.

To prove the surjectivity let $(K,e)$ be a local pair on~$P$. Denote 
by~$E$ the subring of the ring $\End(P)$ generated by~$K$. Since the group~$K$
is abelian, the ring $E$ is commutative. Let us show that the group
homomorphism
\qtnl{l240703a}
E^+\to P,\quad T\mapsto T(e)
\eqtn
is in fact an isomorphism. Suppose that $T(e)=0$ (here and below we treat
$P$ as an additive group). Then
$$
T(e^k)=(Tk)(e)=(kT)(e)=T(e)^k=0
$$
for all $k\in K$. On the other hand, from the definition of a local pair
it follows that each element of $P\setminus O$ where $O=e^K$, is the
difference (in~$P$) of some elements of~$O$. Thus $T(x)=0$ for all $x\in P$,
i.e. $T=0$. So (\ref{l240703a}) is a monomorphism. Next, obviously the
image of~$K$ with respect to~(\ref{l240703a}) equals~$O$. So (\ref{l240703a}) 
is an epimorphism because $O$ generates~$P$ (see above).

Isomorphism~(\ref{l240703a}) takes $K$ to $O$ and 
$E\setminus K$ to~$P\setminus O$. Since $P\setminus O$ is a subgroup
of~$P$, $E\setminus K$ is a subgroup of the group~$E^+$.
This implies that $E\setminus K$ is an ideal of~$E$
($E\setminus K$ is stable with respect
to multiplication by~$K$ and hence by~$E$). Since $K\subset E^\times$,
we conclude that $E\setminus K$ is the only maximal ideal of~$E$. Thus
$E$ is a local ring with $\rad(E)=E\setminus K$ and $E^\times=K$.
Now using isomorphism~(\ref{l240703a}) we come to a local ring $R$ on~$P$
such that $\rad(R)=P\setminus O$, $R^\times=O$ and $1_R=e$. By the
definition of~$R$ we have $T(x)=xT(e)$ for all $T\in E$ and
$x\in P$. It follows that $K=K_R$, whence $(K,e)=(K_R,1_R)$.
Thus the surjectivity of~(\ref{l240703}) and the first part of the
theorem are proved.

To complete the proof we observe first that any isomorphism of two local
commutative rings on~$P$ induces an inner automorphism of the group $\aut(P)$
taking the multiplications by the units of the first ring to those
of the second one. Conversely, let $R_1$ and $R_2$ be local commutative
rings on $P$ such that $\sigma^{-1}K_{R_1}\sigma=K_{R_2}$ for
some $\sigma\in\aut(P)$. Since the multiplicative group of any ring 
transitively acts by multiplications on itself
we assume without loss of generality that
$(1_{R_1})^\sigma=1_{R_2}$. Obviously, the automorphism
$\sigma$ induces a bijection $r\mapsto r_\sigma$ from $R_1^\times$ onto
$R_2^\times$ such that
\qtnl{l250703}
x^{k_1(r)\sigma}=x^{\sigma k_2(r_\sigma)},\quad x\in P,
\eqtn
for all $r\in R_1^\times$ where $k_i(r')$ denotes the automorphism of~$P$
induced by the multiplication by $r'\in R_i^\times$ in the ring~$R_i$
($i=1,2$). Taking $x=1_{R_1}$ we obtain that $r_\sigma=r^\sigma$ for
all $r\in R_1^\times$. So (\ref{l250703}) implies that
$$
(x\cdot r)^\sigma=x^\sigma\cdot r^\sigma,\quad x\in R_1,\ r\in R_1^\times,
$$
where the multiplication on the right-hand (resp. left-hand) side is meant 
in~$R_1$ (resp. in~$R_2$). Since 
$R_1\setminus R_1^\times\subset R_1^\times-R_1^\times$, the last
identity is true for all $r\in R_1$. Thus $\sigma$ is a ring isomorphism
from $R_1$ onto $R_2$ and we are done.\bull

\crllrl{l110703b}
Let $P$ be a finite abelian $p$-group and $K\le\aut(P)$ the first component
of some local pair on~$P$ (i.e. $K$ transitively
acts on the set $P\setminus P_0$ for some group $P_0<P$). Then
\nmrt
\tm{1} $\orb(K_0,P\setminus P_0)=\{xP_0:\ x\in P\setminus P_0\}$ where 
$K_0$ is the Sylow $p$-subgroup of~$K$.
\tm{2} each orbit of $K$ on $\wh P\setminus\{1\}$ is the set difference
of some $K$-invariant subgroup of~$\wh P$ and its maximal $K$-invariant 
subgroup; moreover $P_0^\bot$ is the smallest nontrivial $K$-invariant 
subgroup of~$\wh P$.
\enmrt
\ecrllr
\proof By Theorem~\ref{thrm4} one can assume that $P=R^+$, $K=K_R$
and $P_0=\rad(R)$ for some commutative local ring~$R$.
Then $K_0$ is induced by the action of $1+\rad(R)$ on~$R^+$
(see~\cite[p.355]{MD}), whence statement~(1) follows. Next, 
for $r\in R$ and
$\chi\in \wh{R^+}$ set $\chi^r(x)=\chi(rx)$, $x\in R^+$. Obviously,
$\chi^r\in\wh{R^+}$ and $(\chi^r)^s=\chi^{rs}$ for all $r,s\in R$.
Since $R\setminus\rad(R)=R^\times$ we have 
$\chi^R=\chi^{R^\times}\cup \chi^{\rad(R)}$ for
all $\chi\in\wh{R^+}$. If the character $\chi$ is not principal, then 
$\chi^{R^\times}\cap \chi^{\rad(R)}=\emptyset$.
(Otherwise, $\chi((u-r)x)=0$ for some $u\in R^\times$, $r\in\rad(R)$
and all $x\in R$. Since $u-v\in R^\times$, this implies that
$\chi$ is principal.) Thus any orbit of $R^\times$
on $\wh{R^+}\setminus\{0\}$
is of the form $\chi^{R^\times}=\chi^R\setminus \chi^{\rad(R)}$. Since, 
obviously, $\chi^R$ and $\chi^{\rad(R)}$ are $R^\times$-invariant subgroups 
of $\wh{R^+}$ and $\chi^{R^\times}=\chi^{K_R}$, statement~(2) follows.\bull

\eSuSS{\bf Proof of Theorem~\ref{thrm2}.}
In the conditions of Theorem~\ref{thrm2} let $\A$ be a $K$-primitive
$K$-invariant S-ring over the group~$G$ and $P_0$ a $K$-invariant
subgroup of~$P$ for which $P\setminus P_0$ is an orbit of $K\cap\aut(P)$.
Obviously, $P_0<P$. We consider three cases.
\vspace{2mm}

{\bf Case 1:} $G=P$ and $P_0\ne\{1\}$. The Sylow $p$-subgroup $K_0$ of~$K$
acts on the nontrivial $p$-group $P_0$ as an automorphism group. So there 
exists $x\in P_0\setminus\{1\}$ fixed by~$K_0$. We claim that the group
$H=P_0$ separates the basic set $X=[x]$ of~$\A$ (we observe that $H$ is a
proper $K$-invariant subgroup of~$G$). Indeed, $X\cap H\ne\emptyset$
by obvious reason, $X\setminus H\ne\emptyset$ by the first implication of
Lemma~\ref{t100703a} for $H'=\lg X\rg$, and $H\le\rad(X\setminus H)$
by statement~(1) of Corollary~\ref{l110703b} 
because $X^{K_0}=X$ due to the $K_0$-invariance of~$\A$ and the choice 
of~$X$. Thus $H\in\H_{sep}(\A)$. Since $H$ is $K$-invariant, we are done
by Corollary~\ref{l300703}.
\vspace{2mm}

{\bf Case 2:} $G=P$ and $P_0=\{1\}$. Let $e\in P\setminus\{1\}$.
From the hypothesis of the theorem it follows that $(K,e)$ is a local pair
on~$P$ such that $P\setminus e^K=\{1\}$.  So by Theorem~\ref{thrm4}
without loss of generality we can assume that $G=\F^+$ and
$K=K_\F$ where $\F$ is a field. Set $L$ to be the setwise stabilizer 
of the basic
set $[1_\F]$ of the S-ring~$\A$ in the group~$K$. Since
$\A$ is $K$-invariant, we have $[r]=r[1_\F]$ for all
$r\in\F^\times$. So every basic set of $\A$ is $L$-invariant.
Moreover, if $x,y\in[r]$ where $r\in\F^\times$, then
$$
[1_\F]x^{-1}y=[r]r^{-1}x^{-1}y=[x]x^{-1}yr^{-1}=[y]r^{-1}
=[r]r^{-1}=[1_\F],
$$
whence $x^{-1}y\in L$. Thus, $\S(\A)=\orb(L,\F)$ and we are done.
\vspace{2mm}

{\bf Case 3:} $G\ne P$. Let $\wh\A$ be the S-ring (over $\wh G$)
dual to~$\A$. Since $\rk(\wh\A)=\rk(\A)$, it suffices to prove
that $\rk(\wh\A)=2$. To do this we observe that by Theorem~\ref{t110703a}
the ring $\wh\A$ is both $K$-invariant and $K$-primitive. Set
$H=(P_0Q)^\bot$ with $Q$ the product of all Sylow $q$-subgroups 
of $G$ for $q\ne p$. Then $H$ is a proper $K$-invariant subgroup
of $\wh G$. By Corollary~\ref{l300703} it suffices to verify that
$H\in\H_{sep}(\wh\A)$. To do this we will show that the group~$H$ separates
the basic set $X=[x]$ of~$\wh\A$ where $x\in\wh Q\setminus\{1\}$
(such an~$x$ does exist because $Q\ne\{1\}$). Indeed, obviously
$X\setminus H\ne\emptyset$. So it remains to prove that
\qtnl{l150703b}
X\cap H\ne\emptyset,\quad H\le\rad(X\setminus H).
\eqtn
First we observe that the first formula of~(\ref{l150703b}) is the
consequence of the second one. Indeed, otherwise 
$X=X\setminus H$ and $\rad(X)\ge H$. However, $H$ is 
a proper $K$-invariant subgroup of $\wh G$. So by the second
implication of Lemma~\ref{t100703a} with $G$ replaced by $\wh G$ 
and $H'=\lg X\rg$ we obtain that $\rad(X)=\wh G$ which is impossible.
Let us prove the second formula of~(\ref{l150703b}). To do this
set $L=K\cap\aut(P)$. Then from statement (2) of Corollary~\ref{l110703b}
with $K$ replaced by $L$, it follows that
\qtnl{l150703a}
H\setminus\{1\}\in\orb(L,\wh P)\quad
\textstyle{\rm and}\quad
\rad(O)\ge H\quad
\textstyle{\rm for\ all}\quad
O\in\orb(L,\wh P\setminus H).
\eqtn
Moreover, since the ring $\wh\A$ is $L$-invariant, the set $X$ is
$L$-invariant by the choice of $x$, and so
the set $y\wh P\cap X$ is also $L$-invariant for all $y\in\wh G$.
This implies that for each $y\in\wh Q$ we have
\qtnl{l150703}
(y\wh P\cap X)\setminus H=y\cdot\bigcup_{O\in\O_y}O
\eqtn
for some $\O_y\subset\orb(L,\wh P)$. From (\ref{l150703a}) it follows
that the cardinality of any element of $\orb(L,\wh P)$ other than $\{1\}$
and $H\setminus\{1\}$ is divided by~$p$. By statement~(2) of 
Lemma~\ref{t100703b} this implies that $\{1\}\in\O_y$ iff 
$H\setminus\{1\}\in\O_y$. So using (\ref{l150703a}) once more we conclude
that the radical of the left-hand side of (\ref{l150703}) contains~$H$.
Thus $\rad(X\setminus H)\ge H$ and we are done.\bull

\end{document}